\documentclass{article}
\usepackage{spconf,amsmath,amssymb,graphicx}
\usepackage{bmpsize}
\usepackage{algpseudocode,algorithm}
\usepackage{cite}
\usepackage{color}
\usepackage{comment}
\usepackage{mathtools}
\usepackage{setspace}
\usepackage{url}
\usepackage{booktabs}

\algnewcommand\algorithmicinput{\textbf{Input:}}
\algnewcommand\Input{\item[\algorithmicinput]}
\algnewcommand\algorithmicoutput{\textbf{Output:}}
\algnewcommand\Output{\item[\algorithmicoutput]}

\title{Numerical Spectrum Linking: Identification of Governing PDE \\ via Koopman-Chebyshev Approximation}

\twoauthors
  {Phonepaserth SISAYKEO
  }
	{Graduate School of Sci. \& Tech.\\
  Niigata Univ., Japan
  }
  {Shogo MURAMATSU\sthanks{This work was supported by JSPS KAKENHI Grant Numbers JP22H00512, JP24H00365 and JP24K21314.}}
	{Faculty of Engineering \\
  Niigata Univ., Japan}

\begin{document}
\ninept

\maketitle

\begin{abstract}
A numerical framework is proposed for identifying partial differential equations (PDEs) governing dynamical systems directly from their observation data using Chebyshev polynomial approximation. In contrast to data-driven approaches such as dynamic mode decomposition (DMD), which approximate the Koopman operator without a clear connection to differential operators, the proposed method constructs finite-dimensional Koopman matrices by projecting the dynamics onto a Chebyshev basis, thereby capturing both differential and nonlinear terms. This establishes a numerical link between the Koopman and differential operators. Numerical experiments on benchmark dynamical systems confirm the accuracy and efficiency of the approach, underscoring its potential for interpretable operator learning. The framework also lays a foundation for future integration with symbolic regression, enabling the construction of explicit mathematical models directly from data.
\end{abstract}


\section{Introduction}

According to Shannon's sampling theorem, a band-limited continuous signal can be reconstructed if sampled at a frequency at least twice its maximum frequency\cite{Shannon1949a,Unser2000a,Eldar2014SamplingSystems}. This study applies this concept to the sampling and reconstruction of governing laws (operations) underlying dynamical systems.

As countermeasures against natural disasters such as river flooding, optimized structural design and control of physical phenomena are increasingly required. To realize such countermeasures, models that accurately capture the governing laws of target dynamics are indispensable. However, many natural phenomena are highly complex for example, the relationship between river flow and sediment transport and cannot always be represented by existing formulations~\cite{Moteki2022CaptureBars,Moteki2023OnSandbars,Ohara2024Physics-informedModel}. Furthermore, new discoveries have overturned established laws~\cite{Karisawa2025Acceleration-inducedSlopes}, emphasizing the importance of data-driven identification of governing dynamics for newly observed phenomena.

\begin{figure}[tb]
  \centering
  \includegraphics[width=.9\linewidth]{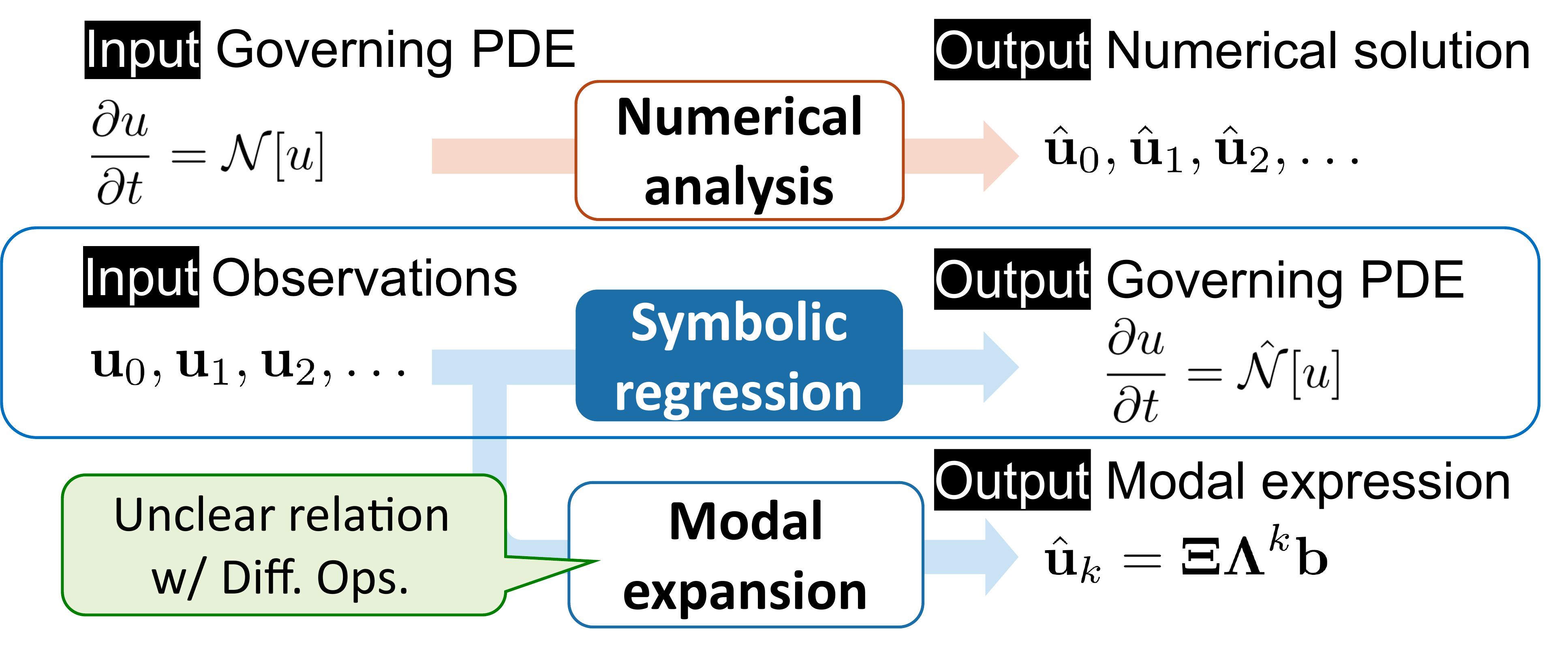}
  \caption{Establishing a foundation for data-driven symbolic regression of PDEs governing dynamics.}\label{fig:background}
\end{figure}

Data-driven modeling has been advanced by mode-expansion methods such as DMD~\cite{Brunton2022Data-DrivenControl} and its variants~\cite{Baddoo2023Physics-informedDecomposition,Kobayashi2023Multi-ResolutionModeling}, which approximate the Koopman operator~\cite{2020TheControl}. While the resulting time-evolution equations are directly applicable to simulation, they lack the interpretability and generality of PDEs, limiting their utility for understanding physical mechanisms and extending to different environments. Therefore, bridging mode expansion with PDE identification, and ultimately symbolic regression~\cite{Wu2025DiscoveringFramework}, is of significant importance as shown in Fig.~\ref{fig:background}.

Pioneering work in PDE identification from observational data includes PDE-FIND by Rudy \textit{et al.}~\cite{Rudy2017Data-drivenEquations}, which employed sparse regression on candidate libraries of differential operators to successfully identify PDEs of complex systems. However, this approach suffers from dependence on predefined libraries and sensitivity to noise. Chen \textit{et al.}'s SGA-PDE~\cite{Chen2022SymbolicSGA-PDE} mitigated library dependence by evolving candidates through genetic algorithms, but at the cost of computational load and stability issues. More recently, Du \textit{et al.}'s DISCOVER~\cite{Du2024DISCOVER:Learning} introduced reinforcement learning to identify PDEs containing unknown terms, but challenges remain regarding computational efficiency, reward design, and robustness to noise.

In contrast, this study focuses on the intrinsic connection between differential operators and the Koopman operator. Conventional DMD has lacked clarity in relating basis functions \(\{\varphi_m\}\) of the observable space to the differential operator \(\mathcal{N}\). By explicitly establishing their numerical connection, this study formulates a problem for PDE identification directly from data. Note that recent advances in surrogate models such as physics-informed neural networks (PINNs) have demonstrated the effectiveness of incorporating known physical laws into loss functions~\cite{Raissi2019Physics-informedEquations,Huang2025PartialSurvey}. However, PINNs presuppose known governing equations, and thus differ fundamentally in purpose from this research.

The objective of this study is to establish a theoretical foundation for numerically identifying governing laws from observational data. Shi \textit{et al.} utilizes the relationship between Chebyshev bases and Koopman operators to enhance the efficiency of numerical computations~ \cite{Shi2024KoopmanStudyb}.
Inspired by the article, a framework termed \textit{numerical spectrum linking} is proposed, combining Chebyshev expansions~\cite{Cohl2013GeneralizationsIntegrals} with Koopman mode expansions.
Chebyshev polynomials approximate observed functions through orthogonal basis expansions, while DCT enables efficient finite-dimensional approximations~\cite{Ahmed1974DiscreteTransform,Ochoa-Dominguez2019DiscreteTransform}, allowing differential operators to be represented as matrix operations on coefficient vectors~\cite{Barrio2004AlgorithmsSeries}.
Meanwhile, Koopman mode expansion describes the time evolution of nonlinear dynamical systems using the linear Koopman operator, which represents complex dynamics within a linear framework. Our challenges lie in connecting these two representations to facilitate PDE identification from data.

The contributions of this study are summarized as follows:
\begin{itemize}
\itemsep 0pt
  \item A theoretical formulation linking Chebyshev polynomial expansions with Koopman eigenfunctions.  
  \item A numerical procedure to connect the physics-informed operator $\mathbf{K}^{\star}$ and the observation-driven operator $\hat{\mathbf{K}}$.  
  \item Validation through experiments on multiple PDEs, showing that the governing dynamics can be identified from limited observational data.  
\end{itemize}

\section{Overview of Foundational Theories}
\label{sec:related}
\color{black}

\begin{figure}[tb]
  \centering
  \includegraphics[width=.98\linewidth]{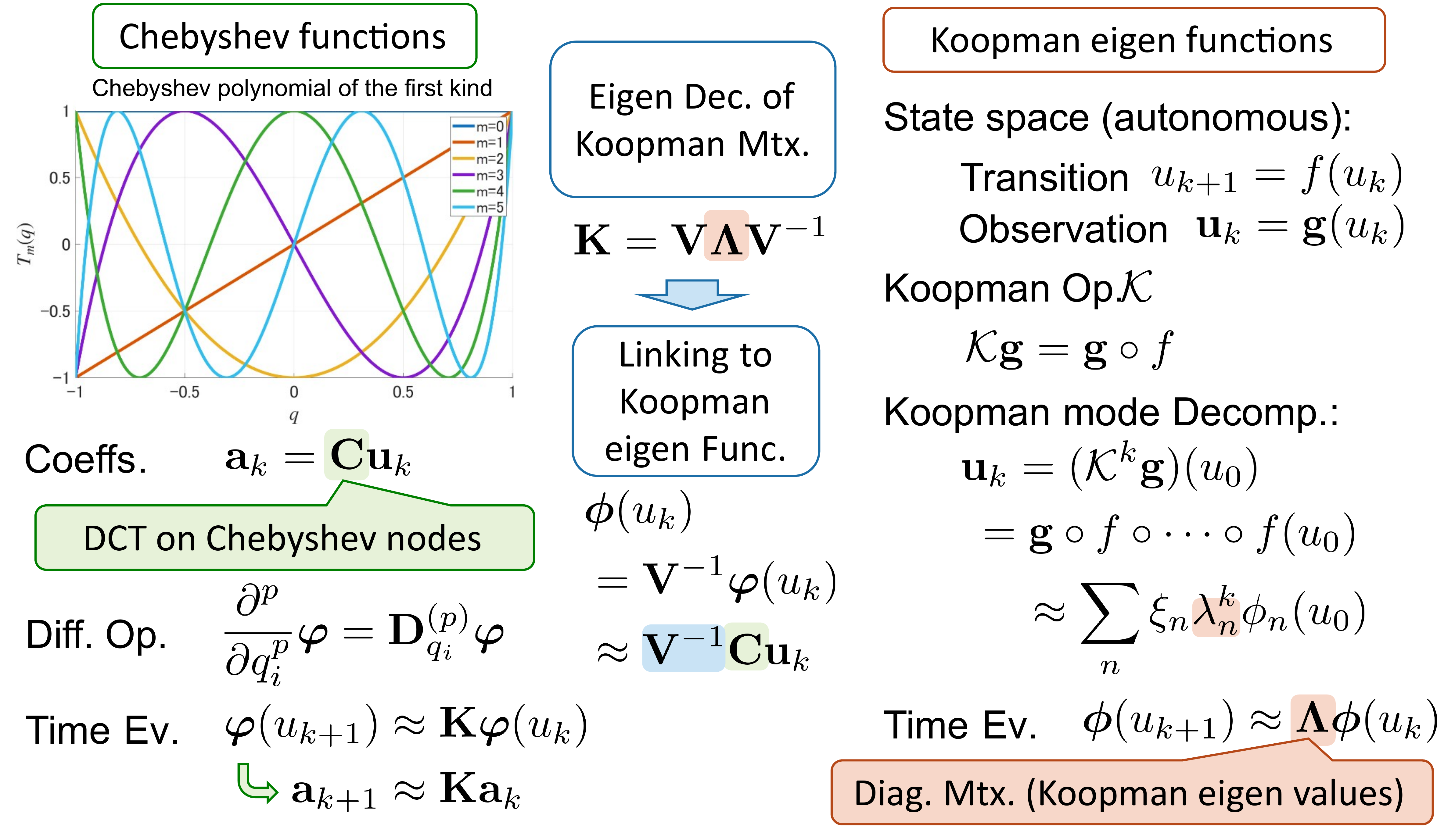}
  \caption{Outline of Chebyshev polynomials and the Koopman operator}\label{fig:outline}
\end{figure}

In this section, we provide foundational theories for the proposed framework. 
We first review the Koopman operator in discrete time and its connection to DMD, which offers a linear representation of nonlinear dynamics. 
We then introduce Chebyshev polynomials and their relation to DCT, which form the basis of efficient spectral approximation. 
Together, these concepts establish the foundation for constructing observation-driven operator models as shown in Fig.~\ref{fig:outline} with improved interpretability and computational efficiency.

\subsection{Koopman Operator}

The description of dynamical systems often benefits from a state-space representation, which is compatible with modern control theory. 
In continuous time, such systems are described by ordinary differential equations (ODEs) or PDEs, 
while in discrete time they are described by difference equations or maps. 
These representations are directly connected to optimal control and state estimation.

A theoretical framework for approximating state-space models is provided by the \emph{Koopman operator}~\cite{2020TheControl,Strasser2025AnGuarantees}. 
The Koopman operator enables one to treat a nonlinear dynamical system as a linear operator acting on observables. 
Consider the nonlinear discrete-time autonomous system
\begin{subequations}
\begin{align}
u_{k+1} &= f(u_k), \quad u_k\in\mathcal{M},\\
{y}_{k} &= {g}(u_{k}), \quad {y}_{k}\in\mathbb{R},
\end{align}
\end{subequations}
where $k\in\mathbb{N}_0$ denotes discrete time, $\mathcal{M}$ is the state space, 
$f:\mathcal{M}\to\mathcal{M}$ is a nonlinear (regular) state-transition map, 
and ${g}:\mathcal{M}\to\mathbb{R}$ is an \textcolor{black}{ observable}.

Let $\mathcal{M}^\ast$ denote the space of admissible observable ${g}$. 
The Koopman operator $\mathcal{K}:\mathcal{M}^\ast\to\mathcal{M}^\ast$ is defined by composition as
\begin{equation}
  \mathcal{K}{g} \;\coloneqq\; {g}\circ f, \quad \forall g\in\mathcal{M}^\ast.    
\end{equation}
Thus, instead of describing the evolution of the state directly, the Koopman operator describes the evolution of observables. 
Note that $\mathcal{K}$ is linear, i.e.,
\(
\mathcal{K}(a {g}_1 + b {g}_2) = a \mathcal{K}{g}_1 + b \mathcal{K}{g}_2,
\)
even though the underlying dynamics $f$ is nonlinear. 
Unless $\mathcal{M}$ is finite, however, the operator $\mathcal{K}$ is infinite-dimensional. 
The operator can also be defined for vector-valued observables $\mathbf{g}:\mathcal{M}\to\mathbb{R}^M$.

\subsubsection{Koopman Mode Expansion}
If the system is integrable, the state space compact, and we consider scalar observables, 
the spectral properties of $\mathcal{K}$ allow for an eigenfunction expansion. 
A function $\phi_\lambda\in\mathcal{M}^\ast\setminus\{0\}$ satisfying
\begin{equation}
  \mathcal{K}\phi_\lambda = \phi_\lambda\circ f = \lambda \phi_\lambda
\end{equation}
is called a \emph{Koopman eigenfunction} with corresponding \emph{Koopman eigenvalue} $\lambda\in\mathbb{C}$. 
Using a family of eigenfunctions $\{\phi_{\lambda_m}\}_{m=0}^\infty$, an observable $g$ can be expanded as
\(
g = \sum_{m=0}^\infty \xi_m \phi_{\lambda_m},
\)
where $\xi_m\in\mathbb{C}$ are the \emph{Koopman modes}. 
The time evolution of $g$ along a trajectory is then
\(
y_{k} = (\mathcal{K}^{k}g)(u_0)
= \sum_{m=0}^\infty \xi_m \lambda_m^{k} \phi_{\lambda_m}(u_0).
\)

For vector-valued observables $\mathbf{g}$, both the measurements 
$\mathbf{y}_{k}\in\mathbb{R}^M$ 
and the Koopman modes $\boldsymbol{\xi}_m\in\mathbb{C}^M$ become vector-valued, yielding
\(
\mathbf{y}_{k} = (\mathcal{K}^{k}\mathbf{g})(x_0) 
= \sum_{m=0}^\infty \boldsymbol{\xi}_m \lambda_m^{k} \phi_{\lambda_m}(
    x_0).
\)

\subsubsection{Koopman Matrix Approximation}
For numerical computation, the infinite-dimensional operator $\mathcal{K}$ is approximated 
in a finite-dimensional subspace spanned by a set of basis functions 
$\{\varphi_m\}_{m=0}^{M-1}\subseteq\mathcal{M}^\ast$. 
Approximating $g(u_k)$ by a linear combination of these basis functions leads to the relation
\begin{equation}
(\mathcal{K}g)(u_k) \approx \sum_{m=0}^{M-1} a_m (\mathcal{K}\varphi_m)(u_k)
= \sum_{m=0}^{M-1} a_m \varphi_m(u_{k+1}). \label{eq:koopmanmatrix0}
\end{equation}
If each $\varphi_m(u_{k+1})$ can be expressed as a linear combination of $\{\varphi_n(u_k)\}$, 
the system admits a matrix representation
\begin{equation}
\boldsymbol{\varphi}(u_{k+1}) \approx \mathbf{K}\,\boldsymbol{\varphi}(u_{k}),
\end{equation}
where $\boldsymbol{\varphi}(u) = (\varphi_0(u),\ldots,\varphi_{M-1}(u))^\top$ and
$\mathbf{K}\in\mathbb{R}^{M\times M}$ is called the \emph{Koopman matrix}.
Given data $\{u_{k}\}_{k=0}^{N-1}$,
$\mathbf{K}$ can be estimated by solving the least-squares problem
\begin{equation}
\hat{\mathbf{K}} = \arg\min_{\mathbf{K}} \sum_{k=0}^{N-2} 
\left\| \boldsymbol{\varphi}(u_{k+1}) - \mathbf{K}\boldsymbol{\varphi}(u_{k})\right\|^2.
\end{equation}
If the basis functions are Koopman eigenfunctions, the matrix $\mathbf{K}$ becomes diagonal. 
Eigenvalue analysis of the estimated $\mathbf{K}$ provides access to Koopman modes, eigenvalues, and eigenfunctions. 
This procedure underlies DMD and its extensions~\cite{Brunton2022Data-DrivenControl,Bistrian2025ReducedLearning}.

\subsection{Chebyshev Polynomials}
In this subsection, we summarize the basic properties of Chebyshev polynomials and their role in spectral approximation. 
Chebyshev polynomials form an orthogonal basis on $[-1,1]$ with respect to a specific weight function, and their tensor-product extension enables efficient approximation of multivariate functions. 
By sampling at Chebyshev nodes and exploiting the correspondence with the type-II DCT (DCT-II), the expansion coefficients can be computed efficiently
~\cite{Ochoa-Dominguez2019DiscreteTransform}. 
This connection provides a practical tool for constructing finite-dimensional operator representations and serves as a foundation for the Koopman operator analysis developed later.

\subsubsection{Multidimensional Chebyshev basis and expansion}
We consider a real-valued function on a $D$-dimensional hypercube,
\begin{equation}
u(\mathbf{q}) \in \mathbb{R}, 
  \quad \mathbf{q} \in [-1,1]^D .
\end{equation}
For simplicity, the time index is omitted.

Let $T_m(q)\coloneqq \cos(m\arccos q)$ be the Chebyshev polynomials of the first kind, $m\in\mathbb{N}_0$, $q\in[-1,1]$.
We define the $D$-dimensional tensor-product basis
\(
  \tau_{\mathbf{m}}(\mathbf{q}) \;=\; \prod_{d=1}^D T_{m_d}(q_d)
  \), 
  \(
  \mathbf{m}=(m_1,\dots,m_D)^\intercal\in\mathbb{N}_0^D.
\)
Then, $u$ admits the (formal) Chebyshev expansion
\(
  u(\mathbf{q}) \;=\; \sum_{\mathbf{m}\in\mathbb{N}_0^D} 
  a_\mathbf{m}\,\tau_{\mathbf{m}}(\mathbf{q}).
\)

\subsubsection{Sampling at Chebyshev nodes and DCT-II computation}
For a finite-dimensional approximation, we select 
\(
\mathcal{N}(\mathbf{M}) = \{0, 1, \dots, M_d-1\}^D
\),
\( 
\mathbf{M}=\mathrm{diag}(M_1,M_2,\dots,M_D),
\)
and write
\begin{equation}
u(\mathbf{q}) \approx \sum_{\mathbf{m}\in\mathcal{N}(\mathbf{M})} a_\mathbf{m}\,\tau_{\mathbf{m}}(\mathbf{q})
= \sum_{\mathbf{m}\in\mathcal{N}(\mathbf{M})} \check{a}_\mathbf{m}\,\check{\tau}_{\mathbf{m}}(\mathbf{q}).
\end{equation}
Here, we rescale with
\(
\check{a}_\mathbf{m} = \gamma_{\mathbf{m}}^{-1} a_\mathbf{m}, \quad
\check{\tau}_{\mathbf{m}} = \gamma_{\mathbf{m}} \tau_{\mathbf{m}}, \quad
\gamma_{\mathbf{m}}=\prod_{d=1}^D \gamma^{(d)}_{m_d},
\)
and choose per-dimension factors 
$\gamma^{(d)}_0 = M_d^{-1/2}$ and $\gamma^{(d)}_{m\ge 1} = (2/M_d)^{1/2}$,
so that the discrete transform below becomes orthonormal.
With this scaling, we have
\(
\check{a}_\mathbf{m}
\;=\;
{\langle \check{\tau}_{\mathbf{m}}, u \rangle_w}/{\langle \check{\tau}_{\mathbf{m}}, \check{\tau}_{\mathbf{m}} \rangle_w}
\;=\;
(\prod_{d=1}^D \frac{M_d}{\pi})\,\langle \check{\tau}_{\mathbf{m}}, u \rangle_w
\),
where $w(q)=\prod_{d=1}^D (1-q_d^2)^{-1/2}$ is the Chebyshev weight.

Assume each $M_d$ is even and use interior Chebyshev nodes (no endpoints)
\(
  p_{n_d} \;=\; \cos({(2n_d+1)\pi}/({2M_d})), 
  \quad n_d=0,\dots,M_d-1 .
\)
Let $\mathbf{p}_{\mathbf{n}}=(p_{n_1},\dots,p_{n_D})^\intercal$ and sample
\(
  \mathbf{u} \;\coloneqq\; \bigl(u(\mathbf{p}_{\mathbf{n}})\bigr)_{\mathbf{n}\in\mathcal{N}(\mathbf{M})}.
\)
Then, the (vectorized) coefficient vector 
$\mathbf{a}=\bigl(\check{a}_\mathbf{m}\bigr)_{\mathbf{m}\in\mathcal{N}(\mathbf{M})}$ 
is well-approximated by a $D$-dimensional DCT-II:
\begin{equation}
  \mathbf{a} \;\approx\; \mathbf{C}\mathbf{u},\label{eq:a_approx_Cu}
\end{equation}
where $\mathbf{C}$ denotes the (Kronecker-structured) matrix of the $D$-dimensional DCT-II.
We can represent the basis map as
\(
\boldsymbol{\varphi}(u) 
\;\approx\; \mathbf{C} \mathbf{u},
\)
which sends samples on Chebyshev nodes to (scaled) Chebyshev coefficients. 
If $\mathbf{u}_{k}$ collects the samples of the state at time $k$,
i.e., $\mathbf{u}_k=\left(u({\mathbf{p}_n}, k\Delta t)\right)_n$, and then the observation relation can be written as
\begin{equation}
  \mathbf{u}_{k} = \mathbf{g}(u_{k}) \approx 
  \mathbf{C}^\intercal\boldsymbol{\varphi}(u_{k}).
\end{equation}

\subsection{Differentiation of Chebyshev polynomials}
We summarize differentiation identities and their implementation in both coefficient space and node space, consistent with the Chebyshev/DCT-II scaling used above.

For the $D$-dimensional tensor-product basis $\tau_{\mathbf{m}}(\mathbf{q})$ with 
$\mathbf{q}=(q_d)_{d=1}^{D}$, the partial derivative with respect to $q_d$ acts on the coefficient vector as a Kronecker product operator
\(
\mathbf{D}^{(d)} 
\;=\; \mathbf{I}_{M_D}\otimes\cdots\otimes \mathbf{I}_{M_{d+1}}
\otimes \mathbf{D}_{M_d}
\otimes \mathbf{I}_{M_{d-1}}\otimes\cdots\otimes \mathbf{I}_{M_1},
\)
so that for the stacked coefficients ${\mathbf{a}}\in\mathbb{R}^{M_1\cdots M_D}$\cite{Bedratyuk2022DerivationsPolynomials},
\begin{equation}
\frac{\partial }{\partial q_d}u(\mathbf{q})\approx\sum_{\mathbf{m}\in\mathcal{N}(\mathbf{M})} \check{a}_\mathbf{m}'\check{\tau
}_\mathbf{m}(\mathbf{q})
\quad\Longleftrightarrow\quad 
\mathbf{a}'\;\approx\;\mathbf{D}^{(d)}\,\mathbf{a},\label{eq:mtxrep_diff}
\end{equation}
where \textcolor{black}{$\mathbf{a}\approx\mathbf{C}\mathbf{u}$ and $\mathbf{a}'\approx\mathbf{C} (\partial\mathbf{u}/\partial q_d)$}.
Then, the partial derivative evaluated at nodes is computed by
\begin{equation}
\frac{\partial}{\partial q_d}\,\mathbf{u}
\;\approx\;
\mathbf{C}^\top\,\mathbf{D}^{(d)}\,\mathbf{C}\,\mathbf{u}.
\end{equation}

Second and higher derivatives are obtained by repeated application in coefficient space, e.g.,
\( 
\frac{\partial^2 }{\partial q_d \partial q_r}u(\mathbf{q})
\approx\sum_{\mathbf{m}\in\mathcal{N}(\mathbf{M})} \check{a}_\mathbf{m}''\check{\tau
}_\mathbf{m}(\mathbf{q})
\)
\(\Longleftrightarrow\) 
\(
\mathbf{a}''\approx\mathbf{D}^{(d)}\mathbf{D}^{(r)}\mathbf{a},
\)
which in node space corresponds to 
$\mathbf{C}^\top \mathbf{D}^{(d)}\mathbf{D}^{(r)} \mathbf{C}$.

\section{Proposed Numerical Spectrum Linking}
\label{sec:method}

\begin{figure}[tb]
  \centering
  \includegraphics[width=.98\linewidth]{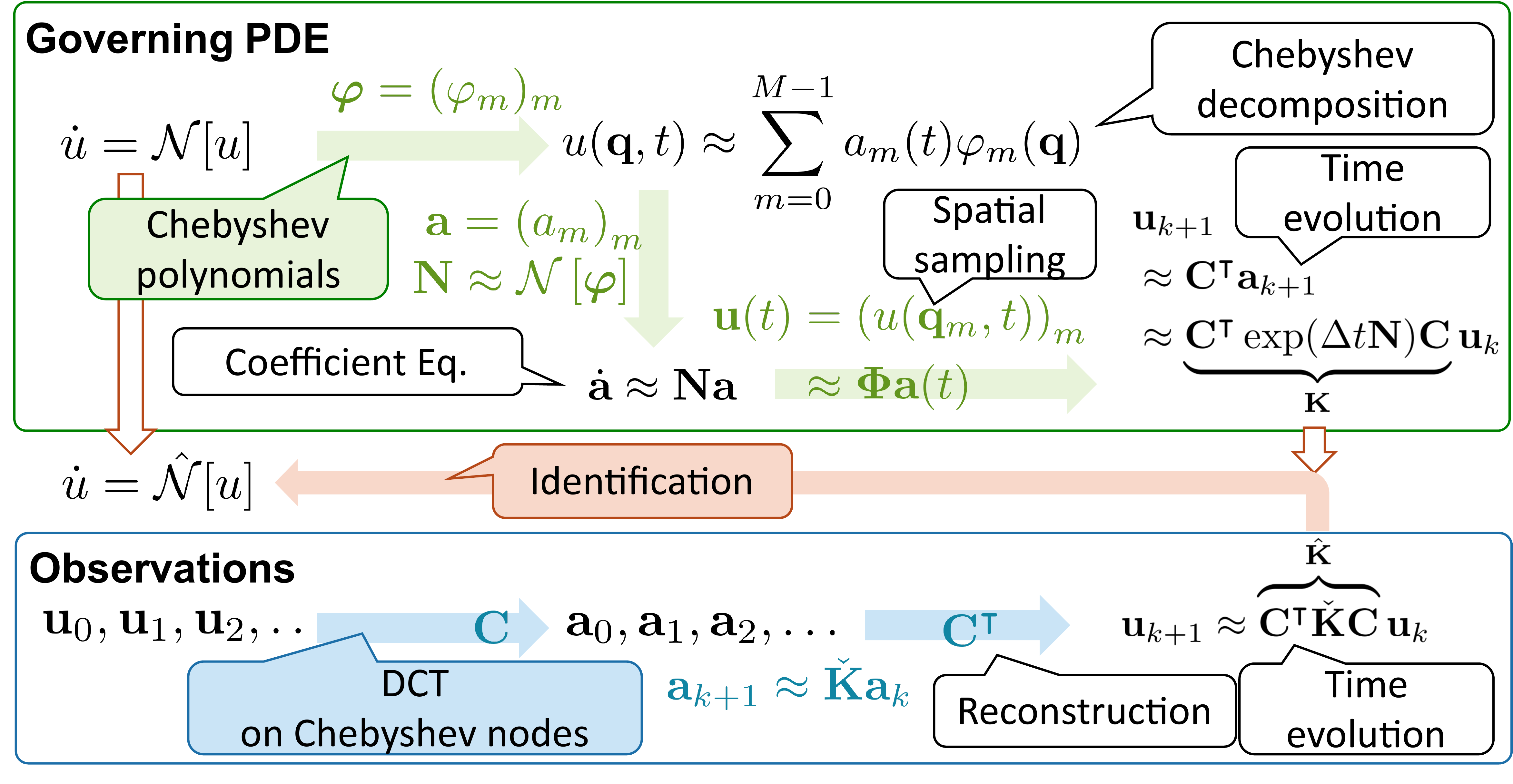}
  \caption{Outline of the proposed numerical spectrum linking method}\label{fig:proposal}
\end{figure}

We present a framework for deriving the Koopman matrix directly from governing equations via Chebyshev approximation, and for comparing it with an operator identified from observed data. The goal is to establish consistency between physics-informed and observation-driven representations of spatiotemporal dynamics.

Consider a governing equation of the form
\begin{equation}
    \dot{u} = \mathcal{N}[u],\label{eq:governingeq}
\end{equation}
where $\dot{u}\coloneqq \partial u/\partial t$, and $\mathcal{N}[\cdot]$ is a spatial differential operator. 

\subsection{Observation-driven identification of Koopman matrix} 

Let us assume that we have observations $\{\mathbf{u}_k\}_{k=0}^{N-1}$ from a dynamics governed by \eqref{eq:governingeq}, where $\mathbf{u}_k=\mathbf{g}(u_k)$. 
Since $\mathbf{a}_{k}\approx\mathbf{C}\mathbf{u}_k\approx\boldsymbol{\varphi}({u}_{k})$ from \eqref{eq:a_approx_Cu}
at time $k$, and the orthonormality of $\mathbf{C}$, the Koopman matrix can be estimated by 
\begin{equation}
    \hat{\mathbf{K}}
    = \arg\min_{\mathbf{K}} 
    \sum_{k=0}^{N-2} 
  \bigl\| \mathbf{a}_{k+1} - \mathbf{K}\,\mathbf{a}_{k}
  \bigr\|_2^2.
\end{equation}
Equivalently, we have
\(
{\hat{\mathbf{K}}}
=\mathbf{A}_1\mathbf{A}_0^{\dagger},
\)
where
\(
\mathbf{A}_i
\coloneqq\bigl(\mathbf{a}_i\quad\mathbf{a}_{i+1}
\)
\(
\ldots\quad\mathbf{a}_{i+N-2}\bigr),
\)
\(
i\in\{0,1\}.
\)

\subsection{Equation-driven derivation of Koopman matrix}

From \eqref{eq:mtxrep_diff}, we may represent the governing equation in \eqref{eq:governingeq} as 
\begin{equation}
\dot{u}\approx
\mathcal{N}\left[
\sum_{\mathbf{m}\in\mathcal{N}(\mathbf{M})} \check{a}_\mathbf{m}\check{\tau
}_\mathbf{m}(\mathbf{q})\right]
\quad\Longleftrightarrow\quad 
\dot{\mathbf{a}}(t)\;\approx\;\mathbf{N}\mathbf{a}(t),
\end{equation}
where $\mathbf{N}$ is the matrix representation of the spatial differential operator $\mathcal{N}$ with respect to the Chebyshev basis.

Approximately, in the discrete form, we have
\begin{equation}
  \mathbf{a}_{k+1} \;\approx\; \mathbf{K}^{\star}\mathbf{a}_k,
\end{equation}
where 
\(
\mathbf{K}^{\star} = \exp(\Delta t \mathbf{N}).
\)

Hence, the Koopman matrix can be numerically derived from the governing differential operator.

\subsection{Numerical Spectrum Linking}

We have obtained two approximations of the Koopman matrix: the data-driven one, denoted by 
\(
  \hat{\mathbf{K}},
\)
and the equation-driven one derived from the governing PDE, denoted by 
\(
  \mathbf{K}^{\star}.
\)
Although they are constructed from different principles, their spectral properties should be consistent.

Let 
\begin{equation}
  \hat{\mathbf{K}}
  \hat{\mathbf{v}}_j 
  = \hat{\lambda}_j
  \hat{\mathbf{v}}_j, 
\end{equation}
\begin{equation}
  \mathbf{K}^{\star} \mathbf{v}_j^{\star} 
  = \lambda_j^{\star} \mathbf{v}_j^{\star}
\end{equation}
be the eigenvalue decompositions of the two matrices. Here, $\lambda_j$ are the Koopman eigenvalues and $\mathbf{v}_j$ are the corresponding Koopman modes in the DCT domain.

The comparison of spectra,
\(
  \{
  \hat{\lambda}_j
  \}
  \quad \text{and} \quad
  \{\lambda_j^{\star}\},
\)
verifies whether temporal growth rates and oscillatory behaviors are consistently captured. In addition, the comparison of Koopman modes,
\begin{equation}
 \hat{\boldsymbol{\xi}}_j
 =\mathbf{C}^\intercal
 \hat{\mathbf{v}}_j
  \;\;\longleftrightarrow\;\; 
  \boldsymbol{\xi}_j^{\star}=
  \mathbf{C}^\intercal\mathbf{v}_j^{\star},
\end{equation}
provides a spatial correspondence between data-driven patterns and PDE-derived structures.

Agreement in both eigenvalues and eigenvectors indicates that 
\(\hat{\mathbf{K}}
\) not only reproduces the correct dynamics but also captures the spatial organization of modes encoded in \(\mathbf{K}^{\star}\).
Thus, our proposed numerical spectrum linking offers a unified framework to validate Koopman matrices through simultaneous comparison of eigenvalues and modes.

\section{Performance Evaluation}
\label{sec:eval}

We evaluate whether the proposed numerical spectrum linking framework can correctly identify the governing PDE from observational data. The comparison is performed between the physics-informed operator $\mathbf{K}^{\star}$ and the observation-driven operator $\hat{\mathbf{K}}
$, focusing on the consistency of their spectral properties.

\subsection{Experimental Conditions}
The experimental setup includes four representative PDEs: \emph{Advection-X}, \emph{Advection-Y}, \emph{Diffusion}, and \emph{Advection–Diffusion}.
These cases form a diverse first-order testbed for assessing the framework's discriminative capability.
All simulations are performed with Chebyshev resolution $M=8$, time step $\Delta t = 5 \times 10^{-4}$, and a total horizon of $T=0.5$. 
Reference solutions are obtained using high-accuracy RK4 integration, while candidate Koopman operators $\mathbf{K}^{\star}$ are derived via Chebyshev spectral discretization of the governing operators. 
Observation-based operators $\hat{\mathbf{K}}
$ are identified by regression from snapshot data. 
The evaluation workflow is summarized in Fig.~\ref{fig:proposal}.

To quantify similarity, two measures are adopted. We define distance and similarity of Koopman matrices as
\begin{equation}
 d
 \big(\mathbf{K}^{\star},\hat{\mathbf{K}}
 \big)
 \coloneqq
 \frac{1}{M}\sum_{i=1}^M \min_j 
 \sqrt{
 \big\|\lambda_i^{\star}\mathbf{v}_i^{\star} - \hat{\lambda}_j
 \hat{\mathbf{v}}_j
 \big\|_2^2
 },
 \end{equation}

\begin{equation}
s
\big(\mathbf{K}^{\star},\hat{\mathbf{K}}
\big)
\coloneqq \frac{1}{M}\sum_{i=1}^M \max_j 
\frac{\big|
\langle\lambda_i^{\star}\mathbf{v}_i^{\star},\hat{\lambda}_j
\hat{\mathbf{v}}_j
\rangle
\big|} {\|\lambda_i^{\star}\mathbf{v}_i^{\star}\|_2\,\|\hat{\lambda}_j
\hat{\mathbf{v}}_j
\|_2}.
\end{equation}

\subsection{Analysis Results}
To evaluate identifiability, confusion matrices are constructed where each row corresponds to the true PDE used to generate the observation data and each column corresponds to a candidate operator. Diagonal dominance indicates successful identification, since $d$ should be minimized and $s$ maximized when $\hat{\mathbf{K}}$ is compared with its true $\mathbf{K}^{\star}$. 
Table~\ref{tab:eigdist} reports the mean minimal distance $d$, where the smallest value in each row consistently appears on the diagonal, confirming correct identification of the underlying PDE. 
Table~\ref{tab:vecsim} presents the
similarity $s$, showing that diagonal entries are the largest across all cases, indicating strong spatial consistency between observation-driven and physics-informed operators. 

\begin{table}[tb]
  \centering
  \caption{Confusion matrix based on Koopman matrix distance $d$  (lower$\downarrow$ is better). Rows: true PDE, Columns: candidate PDE.
  }
  {\footnotesize
  \label{tab:eigdist}
  \begin{tabular}{lccccc}
    \toprule
    Candidate$\backslash$True & Adv-X & Adv-Y & Diffusion & Adv-Diff \\
    \midrule
    Advection-X        & \textbf{0.93180}    & 0.97184            & \underline{0.97099}& 0.96243 \\
    Advection-Y        & 0.98997            & \textbf{0.90342 }  & 0.98023            & 0.96240 \\
    Diffusion          & \underline{0.95643}& \underline{0.94818}& \textbf{0.95946}   & \underline{0.85459} \\
    Advection-Diffusion&0.96269         & 0.95783            & 0.97507            & \textbf{0.81905} \\
    \bottomrule
  \end{tabular}
  }
\end{table}

\begin{table}[tb]
  \centering
  \caption{Confusion matrix based on Koopman matrix similarity $s$ (higher$\uparrow$ is better). Rows: true PDE, Columns: candidate PDE.
  }
  \label{tab:vecsim}
  {\footnotesize
  \begin{tabular}{lccccc}
    \toprule
    Candidate$\backslash$True & Adv-X & Adv-Y & Diffusion & Adv-Diff \\ 
    \midrule
    Advection-X        & \textbf{0.31580}    & 0.28191           & \underline{0.23152}& 0.25326 \\
    Advection-Y        & 0.16971            & \textbf{0.40853}  & 0.19622            & 0.26286 \\
    Diffusion          & \underline{0.27840} & \underline{0.30830}& \textbf{0.26927}   & \underline{0.42913} \\
    Advection-Diffusion& 0.26172            & 0.27663           & 0.21660             & \textbf{0.55555} \\
    \bottomrule
  \end{tabular}
  }
\end{table}

\section{Conclusion}
\label{sec:conclusion}

This work has presented a numerical spectrum linking framework for identifying governing PDEs by comparing Koopman operators derived from governing equations with those obtained from observational data via Chebyshev approximation. Numerical experiments on advection, diffusion, and advection–diffusion equations demonstrated that the proposed spectral distance and similarity measures yield diagonally dominant confusion matrices, thereby confirming reliable PDE identification even at modest resolution. 

Future research will extend this framework to higher-dimensional problems and investigate robustness under noisy conditions. In addition, the treatment of second-order PDEs such as the wave equation will be addressed by introducing an augmented state representation (e.g., $[u, u_t]$), which lies beyond the scope of the present study but represents an important direction for further development.

\newpage
\bibliographystyle{IEEEtran}


\end{document}